\documentstyle[11pt]{article}
\author{Vyacheslav A. Artamonov, Alexander A. Totok}
\title{Actions of Pointed Hopf Algebras}
\date{}
\newtheorem{theorem}{Theorem}[section]
\newtheorem{lem}[theorem]{Lemma}
\newtheorem{prop}[theorem]{Proposition}
\newtheorem{dfn}[theorem]{Definition}
\newtheorem{que}[theorem]{Question}
\newcommand{\demo}{{\bf Proof. }}
\newcommand{\edemo}{\quad \mbox {$\Box$} }
\begin{document}
\maketitle
\section{Introduction}
Throughout this paper $H$ is a finite dimentional Hopf algebra over 
a field $k$, and $A$ is a associative $k$-algebra. 
\begin{dfn} It is said that $H$ acts on $A$, if $A$ is
left $H$-module and for any $h \in H, \quad a,b \in A$
$$h(ab) = \sum_h (h_{(1)} a)(h_{(2)} b), 
\quad h1 = \varepsilon(h),$$
where $\varepsilon : H \to k$ - counit 
and $\Delta$ - comultiplication:
$$\Delta(h) = \sum_h h_{(1)} \otimes h_{(2)} \in H \otimes H.$$
Often  $A$ is called  $H$-module algebra.
\end{dfn}
\begin{dfn} The invariants of $H$ in $A$ is the set $A^H$ of those 
 $a \in A$, that $ha = \varepsilon(h)a$ for each $h \in H$.
\end{dfn}

Straightforvard computations shows, that $A^H$ 
is the subalgebra of $A$. We refer reader to \cite {Mo}, 
\cite {S} for the basic information concerning Hopf algebras 
and their actions on associative algebras.

As a generalization of the well-known fact for group actions  
the following question raised in \cite{Mo} ( Question 4.2.6.) 
\begin{que} \label{mainque}
If $A$ is a commutative $k$-algebra and $H$ any finite
dimensional Hopf algebra such that $A$ 
is $H$-module algebra, is $A$ integral over $A^H$ ?
\end{que}

If $A$ is an affine algebra, then Artin-Tate 
lemma ensures that $A^H$ is also affine.

Some positive answers to question \ref{mainque} are known.
\begin{theorem}[\cite{FS}]  \label{coco}
Let $H$ be a finite dimensional cocommutative 
Hopf algebra and let $A$ 
be a commutative $H$-module algebra. 
Then $A$ is integral extension of $A^H$.
\end{theorem}

Some results on affine invariants were obtained without using
integrality.
\begin{dfn} Element $t \in H$ is called left integral, if 
$ht=\varepsilon(h)t$ for all $h \in H$.
\end{dfn}

Note that $tA \subseteq A^H.$
\begin{theorem}[\cite{Mo}, Theorem 4.3.7]     \label{trace}
Let $A$ be left Noetherian ring which is an affine $k$-algebra,
and assume that $A$ is an $H$-module algebra, such that $tA=A^H$.
Then $A^H$ is $k$-affine.
\end{theorem}

As mentioned in \cite{Mo}, p. 48, if $H$ is semisimple,
then $tA=A^H$. By Maschke's theorem, $H$ 
is semisimple if and only if
$\varepsilon(t) \ne 0$ for some non-zero left integral $t$. 
Since the space $\int^l_H$ of left 
integrals in $H$ is one-dimensional 
(\cite{Mo}, \S 2.2), the semisimplicity 
of $H$ is equivalent to the fact 
that $\varepsilon(t) \ne 0$ for all non-zero left integrals.

The main result of the work is stated in the theorem \ref{main}. 
It gives positive answer to Question \ref{mainque} in some 
partial cases. Let $H$ be a pointed Hopf algebra and let $A$ 
be an affine $H$-module algebra; if one of three conditions 
is satisfied, then $A$ is integral over $A^H$: 
\begin{enumerate}
\item $H$ is commutative as an algebra,
\item char $k = p > 0$,
\item $A$ is integral domain.
\end{enumerate}

We recall that Hopf algebra $H$ is called {\it pointed} 
if every simple subcoalgebra of $H$ is one-dimentional; 
pointed Hopf algebra is called {\it connected} if it has 
only one simple subcoalgebra (one-dimensional). The examples 
of pointed Hopf algebras are given by group algebras, 
universal enveloping algebras. In fact, if $G$ - group, then the 
only simple subcoalgebras of $kG$ are those of the 
form $kg$, $g \in G$. At the same time universal enveloping 
algebras are examples of connected Hopf algebras: the only 
simple subcoalgebra of the universal enveloping algebra is $k1_H$. 

Another important example of pointed Hopf algebras 
is represented by series of Hopf algebras $A_{N,\xi}$, 
where $N \geq 2$ -- integer number and $\xi$ -- root of unity 
of degree $N$, considered in \cite{Ma} (see also section 3 
of this paper). Note that with 
$N=2, \xi=-1$, char $k \ne 2$ algebra $A_{2,-1}$ 
(sometimes called $H_4$) is the only Hopf algebra of minimal 
dimension neither cocommutative nor commutative 
(\cite{Mo}, example 1.5.6). In other words, algebra $A_{2,-1}$ 
is the minimal Hopf algebra not covered by theorem \ref{coco}, 
but covered by theorem \ref{main}.

Notice that ideal generated by $x$ in $A_{N,\xi}$ 
is  nilpotent, therefore it lies in radical of $A_{N,\xi}$, 
i.e., $A_{N,\xi}$ is not semisimple. 
The same fact can be verified using argument of left integrals
of $H$. Thus, theorem \ref{main} is not covered by theorem 
\ref{trace}. 

In spite of numerous partial positive results it turned 
out that hypothesis \ref{mainque} isn't true in general, e
ven for pointed Hopf algebras. The counterexamples is built 
in section 3 for series of pointed Hopf algebras 
$A_{N,\xi} , N \geq 2$, mentioned at previous paragraphs
\footnote{After the work was done authors became aware 
through e-mail by Susan Montgomery, that counterexample 
to hypothesis \ref{mainque} was independently built by 
Zhu Shenglin. Also he obtained some positive results, 
solving problem \ref{mainque}. His paper "Integrality 
of module algebras over its invariants" should have 
been appeared at J.Algebra in 1996.}
.
\section{The main theorem}
The proof of the theorem is based on properties 
of coradical filtration of an arbitrary coalgebra. 
We recall the basic facts.
\begin{dfn} Coradical $C_0$ of coalgebra $C$ is the (direct) 
sum of all simple subcoalgebras of $C$. Further by the 
induction for each $n \ge 1$ define 
$$C_n = \Delta^{-1} ( C \otimes C_{n-1} + C_0 \otimes C )$$
\end{dfn}
\begin{theorem}[ \cite{Mo}, Theorem 5.2.2] 
$\{ C_n \}_{n \ge 0}$
is the family if subcoalgebras with the following properties:
\begin{enumerate}
\item $C_{n-1} \subseteq C_n, \quad 
C = \bigcup\limits_{n \ge 0} C_n,$
\item $\Delta C_n \subseteq \sum\limits_{i=0}^{n} 
C_i \otimes C_{n-i}$.
\end{enumerate}
\end{theorem}

Reader may find more detailed description 
of coradical filtration in \cite{S}, chapter 9.
 
Let $H$ be a pointed finite dimentional Hopf algebra 
over field $k$. Let $G=G(H)$  denote the set of {\it grouplike} 
elements of $H$, i.e.,
$$G=\{\quad g\in H\setminus 
0 \quad |\quad \Delta g=g\otimes g\quad \}.$$
It is known that elements of $G$ are lineary independent, 
$G$ is the group under multiplication arised from multiplication 
in $H$, and subalgebra generated by $G$ is group Hopf algebra $kG$. 
Also $kG$ is coradical of $H$.

By lemma 5.2.8 \cite{Mo}, coradical filtration \{$H_n$\} 
of Hopf algebra $H$ is Hopf filtration of Hopf algebra, i.e.,
$\Delta H_n \subseteq \sum_{i=0}^n H_i \otimes H_{n-i},$ \quad
$H_nH_m \subseteq H_{n+m},$ \quad $SH_n \subseteq H_n$ for 
all $n,m \ge 0$, if and only if $H_0$ is sub-Hopf algebra of $H$.
If $H$ -- pointed finite dimensional Hopf algebra, 
then this condition is obviously satisfied. Moreover, 
coradical filtration is finite.

By theorem 5.4.1 from \cite{Mo} (see also \cite {TW}, \cite{Mi} 
for reference), coradical filtration of pointed Hopf 
algebra $H$ have additional properties. If $x \in H_m,$  
$m \ge 1$, then
$$x = \sum_{g, h \in G(H)} c_{g,h},$$
where
\begin{equation} \label{filtr}
\Delta(c_{g,h}) = c_{g,h} \otimes g + h \otimes c_{g,h} + w,
\end{equation}
$$w \in H_{m-1} \otimes H_{m-1}.$$
Note that if $a,b,g,h \in G$ and $c=ac_{g,h}b$, 
then by (\ref{filtr})
\begin{eqnarray}  \label{g}
\Delta (c)=c\otimes agb+ahb\otimes c+w',\quad  
w'\in H_{n-1}\otimes H_{n-1}.       
\end{eqnarray}         

Define $H^+ = \ker \varepsilon, \quad H_r^+ = H_r \cap H^+$.
Let $A^G$  denote subalgebra of $G$-invariants 
in $A$ ($A^H \subseteq A^G$). Extension $A/A^G$ is integral 
by Noether's theorem for group actions 
($H$ -- finite dimensional, therefore $G$ -- finite group).

Before we start to prove main theorem we are going to obtain few 
auxiliary results.

Let $I$ denote the ideal in $H$ generated by elements of form 
$g-1$ that $g \in G$.
\begin{prop}   \label{1}
$I$ is Hopf ideal in $H$.
\end{prop}
\demo
If $g \in G$, then
\begin{eqnarray*}
&\Delta (g-1)=g\otimes g - 1\otimes 1=\\  
&(g-1)\otimes g +1\otimes (g-1)\in I\otimes H+H\otimes I. 
\end{eqnarray*}
$S( I ) \subseteq I$, because $S(g-1)=g^{-1}-1$ and $S$ is 
anti-homomorphism. This yields the proposition. \edemo

\begin{prop}   \label{2}
If $J$ -- Hopf ideal in $H$, then $H/J$ -- pointed Hopf algebra.
Moreover, natural epimorphism of Hopf algebras $\pi : H \to H/J$ 
induces epimorphism of groups of grouplike 
elements $G(H) \to G(H/J).$
\end{prop}
\demo
This statement is direct consequence 
of corollary 5.3.5 from \cite{Mo}.
\edemo

\begin{theorem}  \label{3}
Let one of three following conditions be satisfied:
\begin{enumerate}
\item char $k = p > 0$.
\item $A$ is integral domain;
\item $H$ -- connected and commutative;
\end{enumerate}
Then there exists the chain of subalgebras 
$A=A_{-1}\supseteq A_0\supseteq \dots \supseteq A_n$
with following properties:
\begin{enumerate}
\item each extension $A_i\supseteq A_{i+1}$ is integral;
\item if $x \in H_i^+$ then $x(A_i)=0.$
\end{enumerate}
\end{theorem}
\demo 
To construct this chain we start with defining $A_0=A^G.$ 
Both of necessary conditions are satisfied.
Let the chain
$$A=A_{-1}\supseteq A_0\supseteq \dots \supseteq 
A_r, \quad r\geq 0,$$ be already constructed and 
let $x \in H_{r+1}^+.$ By (\ref{filtr}) we may assume 
that $x=c_{g,h},$ where $g,h\in G.$ Then
\begin{equation} \label{del1}
\Delta (x)=x\otimes g + h\otimes x +\sum u_j\otimes v_j,
\end{equation}
where $v_j,u_j \in H_r$. As $x \in H^+$, then by (\ref{del1})
$$x=(1\otimes \varepsilon )\Delta (x)=x+\sum u_j\varepsilon (v_j),$$
$$x=(\varepsilon\otimes 1)\Delta (x)=x+\sum \varepsilon(u_j)v_j.$$
Therefore,
$$\sum \varepsilon(u_j)v_j=\sum u_j\varepsilon (v_j)=0,$$
and finally
$$\sum (u_j-\varepsilon(u_j))\otimes (v_j-\varepsilon (v_j))=
\sum u_j\otimes v_j,$$
i.e., we may assume that $u_j,v_j\in H_r^+.$ 

If char $k = p > 0$, then we define $A_{m+1} = A_m^p$. 
Really, by (\ref{del1}) 
$$x(a^p) = h(a)^{p-1}x(a) + h(a)^{p-2}x(a)g(a) + \dots $$
$$\dots + x(a)g(a)^{p-1} = pa^{p-1}x(a) = 0.$$

Suppose $A$ is integral domain and char $k = 0$.
If $a \in A_r, b \in A$, then by (\ref{del1})
$$x(ab) = x(a)b + ax(b),$$
i.e., $x : A_r \to A$ is derivation. By normalization lemma 
(see \cite{B}, chapter 5, \S 3, p.344 ), there exists 
subalgebra of polynomials $k[T_1, \dots, T_d]$ in $A_r$ such 
that extension $A_r / k[T_1, \dots, T_d]$ is integral. 
In this case we have for each $f \in k[T_1, \dots, T_d]$: 
$$x(f) = \sum_{i=1}^d \frac{\partial f}{\partial T_i} a_i,
\quad a_i \in A.$$
Therefore, for each integer $q\geq 1$
$$x^q(f) = \sum_{i_1 + \dots +i_d = q \atop i_s \geq 0} 
\frac{\partial^q f}{\partial T_1^{i_1} \dots \partial T_d^{i_d}}
a_1^{i_1} \dots a_d^{i_d} \quad +$$ 
\begin{equation}     \label{int1}
+ \quad \sum_{j_1 + \dots + j_d = l < q \atop j_s \geq 0} 
\frac{\partial^l f}
{\partial T_1^{j_1} \dots T_d^{j_d}} \Psi_{j_1, \dots, j_d},
\end{equation}
where $\Psi_{j_1, \dots, j_d}$ is the sum of monomials of form
$$x^{\alpha_{j_1}}(a_{j_1}) \dots x^{\alpha_{j_d}}(a_{j_d}), \quad
\alpha_{j_1} + \dots + \alpha_{j_d} = q - l.$$
$H$ is finite dimensional, that is why there exists 
such integer $m$ that 
\begin{eqnarray}    \label{int2}
x^m = \sum_{i = 1}^{m-1} \beta_i x^i, \quad \beta_i \in k.
\end{eqnarray}
Thus, by (\ref{int1}) and (\ref{int2}) for each $f \in k[T_1, 
\dots, T_d]$ and each $i = 1, \dots, d$ we have:
\begin{equation}     \label{int3}
\frac{\partial^m f}{\partial T_i^m} a_i^m + 
\sum_{1 \leq l < m} \frac{\partial^l f}{\partial T_i^l} \Phi_l + 
\Lambda = 0,
\end{equation}
where $\Lambda$ is the sum of all summands from (\ref{int1}) and 
(\ref{int2}) containing 
$$\frac{\partial^q f}{\partial T_1^{j_1} \dots \partial T_d^{j_d}}$$
as a multiplier, and besides one of the 
coefficients $j_s, s \neq i$, is not equal to zero.
Substituting to (\ref{int3}) successivly
$1, T_i, T_i^2, \dots, T_i^m$, we get that 
$$\Lambda = \Phi_1 = \dots =\Phi_{m-1} = a_i^m = 0.$$
We use here the fact that char $k = 0$. As $A$ -- integral domain, 
so $a_i = 0$. Thus $x(k[T_1, \dots, T_d]) = 0$ and we define 
$A_{r+1} = k[T_1, \dots, T_d]$.   

Suppose $H$ -- connected, commutative Hopf algebra and char $k = 0$, 
i.e., $g = h = 1$. By (\ref{del1}),
\begin{equation}   \label{del4}
\Delta(x)=x\otimes 1+1\otimes x + \sum u_j\otimes v_j,
\end{equation}
where $u_j,v_j\in H_r^+.$ 

We consider ideal $HH_r^+$ in $H$.
\begin{lem} \label{4}
$HH_r^+$ is coideal in $H$.
\end{lem}
\demo
Let $u\in H,v\in H_r^+.$ Since $H_r$ -- coalgebra, $H_r^+$ -- it's 
coideal, i.e.,
$$\Delta (v)\in H_r\otimes H_r^+ + H_r^+\otimes H_r,$$
and therefore,
$$\Delta(uv)=\Delta(u)\Delta(v)\in H\otimes 
(HH_r^+) + (HH_r^+)\otimes H. \quad \Box$$

Assume that $x(A_r)\ne 0$. $HH_r^+$ acts as zero on $A_r$, so
$x\notin HH_r^+$. Suppose  $1,x, \dots, x^{m-1}$ 
are linary independent modulo $HH_r^+$ and 
\begin{eqnarray}   \label{del5}
x^m=\sum_{j=0}^{m-1}\alpha_jx^j+w,\quad w\in HH_r^+.
\end{eqnarray}
Choose $k$-basis $e_1,\dots, e_d$ in $H$ such that first elements
$e_1,\dots, e_{t-1}$ form $k$-basis for $HH_r^+$, $e_t=1$,
and $e_{t+1}=x, \dots, e_{t+m-1}=x^{m-1}, d \geq t+m-1$. 
Use $\Delta$ on (\ref{del5}). By (\ref{del4}), lemma \ref{4} 
and commutativity of $H$, 
\begin{equation}   \label{del6}
(x\otimes 1+1\otimes x)^m=
\sum_{j=0}^{m-1}\alpha_j(x\otimes 1+1\otimes x)^j+w',
\end{equation}
$$w'\in H\otimes HH_r^+ + HH_r^+\otimes H.$$
Note that commutativity of $H$ is necessary only here to show, 
that $w'$ really lies in $H\otimes HH_r^+ + HH_r^+\otimes H$.
For each integer $q\geq 1$
$$(x\otimes 1+1\otimes x)^q=\sum_{i=0}^q{q\choose i} 
x^i\otimes x^{q-i}.$$
Subtracting from (\ref{del6}) equation
$$1\otimes x^m=\sum_{j=0}^{m-1}\alpha_j1\otimes x^j+1\otimes w,$$
by (\ref{del5}), ({\ref{del6}) we get  
\begin{equation} \label{del7}
\sum_{i=1}^m{m\choose i}x^i\otimes x^{m-i}=
\sum_{j=1}^{m-1}\alpha_j\sum_{i=1}^j{j\choose i}
x^i\otimes x^{j-i}+ w'',
\end{equation}
$$w''\in H\otimes (HH_r^+) + (HH_r^+)\otimes H.$$
Since char $k=0, so {m\choose 1}=m \ne 0$ in $k$.
From this and (\ref{del7}) it follows that element
$x\otimes x^{m-1}= 
e_{t+1} \otimes e_{t+m-1}$
in $H\otimes H$ is linear combination
of elements $e_s\otimes e_{s'}$, where either $s$ is less then $t+1$ 
or $s'$ is less then $t+m-1$. But it is impossible, because elements 
$e_q\otimes e_{q'}, \quad q,q'=1,\dots d$, 
form basis of $H\otimes H.$
This contradiction shows that $x(A_r)=0.$ So in this case
we define $A_{r+1}=A_r=\cdots = A^G.$ 
\edemo

Notice that reasoning shown above 
is close to that used in \cite{Mo}, \S 5.5, \S 5.6. 

We have come to main
\begin{theorem} \label{main}
Let $A$ be an affine $H$-module algebra, $H$ -- finite dimensional 
pointed Hopf algebra and one of three conditions is satisfied:
\begin{enumerate}
\item char $k = p > 0$;
\item $H$ -- commutative;
\item $A$ -- integral domain.
\end{enumerate}
Then extension $A/A^H$ is integral. 
\end{theorem}
\demo
Assume that $H$ is commutative, then $A^G$ is $H$-module algebra. 
It is sufficiently to show that $A^G$ is stable under $H$-action.
In fact, for each $x \in H, \quad a \in A$,
$$gx(a) = xg(a) = x(a),$$
i.e., $\quad x(a) \in A^G$.
Let $I$ denote the ideal in $H$ generated by 
the elements of form $g-1$ ($g \in G$). By proposition \ref{1} $I$ 
is Hopf ideal. Obviously, it acts as zero on $A^G$. 
Hopf algebra $H/I$ is pointed by proposition \ref{2}, moreover, 
it is connected. Thus the second case of this theorem is reduced 
to consideration of action of connected, commutative Hopf 
algebra $H/I$ on algebra $A^G$. Now we apply theorem \ref{3} 
to all cases. Let
$$A=A_{-1}\supseteq A_0\supseteq \ldots \supseteq A_n$$
be constructed chain of subalgebras. By condition 1) 
of theorem \ref{3}, extension $A/A_n$ is integral; 
by condition 2) $A_n\subseteq A^H.$ 
\edemo

Note that if char $k = p > 0$, 
then $(A^G)^{p^{\dim H}}\subseteq A^H$.
If $H$ is commutative and char $k = 0$, then $A^H=A^G$.
\section{Counterexample to hypothesis}
{\bf Example 3.1} Hopf algebra $H$ may be any from series
$A_{N,\xi}$, \quad $N \ge 2$. Algebra $A_{N,\xi}$ is generated by 
elements $g, x$ with relations
\begin{eqnarray} \label{del8}
g^N = 1, \quad x^N = 0, \quad xg = \xi{}gx,
\end{eqnarray}
where $\xi \in k$ -- root of unity of degree $N$.
Hopf algebra structure on $A_{N,\xi}$ is given as follows:
$$\Delta(g) = g \otimes g, \quad \varepsilon(g) = 1, \quad
S(g) = g^{N-1} = g^{-1},$$
$$\Delta(x) = g \otimes x + x \otimes 1, \quad \varepsilon(x) = 0,
\quad S(x) = -g^{N-1}x.$$
We demand $\xi \not= 1$ and char $k$ = 0.
Algebra $A_{N,\xi}$ is pointed,
$$G(A_N) = \{ 1,g ,g^2, \dots , g^{N-1} \},$$ 
it is non-commutative and non-cocommutative.

Let $A$ be commutative algebra generated by elements $y, z$
with relation $z^2 = 0$.
Define action $A_{N,\xi}$ on $A$:
$$g(y^n) = y^n, \quad g(y^nz) = \xi^{-1}y^nz , \quad
x(y^n) = ny^{n-1}z, \quad x(y^nz) = 0.$$

Straightforvard computations shows the correctness of this action,
i.e., $A_{N,\xi} (I) \subseteq I$, where I is ideal of free algebra
$k<y,z>$ generated by elements 
$yz-zy, z^2.$
We check that relations (\ref{del8}) in $H$ are satisfied:
$$\xi{}gx(y^n) = \xi\xi^{-1}ny^{n-1}z = ny^{n-1}z = xg(y^n),$$
$$\xi{}gx(y^nz) = 0 =xg(y^nz),$$
$$g^N(y^n) = y^n, \quad g^N(y^nz)  = \xi^{-N}y^nz = y^nz,$$
$$x^2(y^n) = nx(y^{n-1}z) = 0, \quad x^2(y^nz) = x(0) = 0,$$
i.e., $x^N(a) = x^2(a) = 0$ for any $a \in A$.

Obviously $A^G = k[y]$ and $A^H = k[y] {} \bigcap$ ker $x = k$ .
But extension $A/A^H$ is not integral, because $A$ 
is not finite $k$-module ( $\dim_k A = \infty$ ).
\section{Conclusion}
As it was shown in example 3.1, hypothesis \ref{mainque} is not 
true in general. Nevertheless all known examples of 
Hopf algebra action shows that if $A$ -- affine, then $A^H$ 
is also affine algebra, although extension $A/A^H$ is not 
always integral. So we ask
\begin{que}
Finite dimensional Hopf algebra $H$ acts on commutative 
\newline affine algebra $A$.
Is $A^H$ affine?
\end{que}


\begin{thebibliography}{99}
\bibitem{B}
N. Bourbaki, "Elements of Mathematics. Commutative Algebra,"
Hermann \& Addison-Wesley, 1972.
\bibitem{FS}
Walter R. Ferrer Santos, Finite generation of the invariants 
of finite dimensional Hopf algebras, {\it J. Algebra} {\bf 165}, 
No. 3 (1994), 543 -- 549.
\bibitem{Ma}
A. Masuoka, Cleft extensions for a Hopf algebra generated 
by a nearly primitive element, {\it Comm. Algebra} {\bf 22}, 
No. 1 (1994), 4537 -- 4559.
\bibitem{Mi}
A. Milinski, Actions of pointed Hopf algebras on prime algebras,
{\it Comm. Algebra} {\bf 23}, No. 1 (1995), 313 -- 333.
\bibitem{Mo}
S. Montgomery, "Hopf Algebras and Their Actions on Rings,"
CBMS, No. 82, Amer. Math. Soc., 1993.
\bibitem{S}
M. Sweedler, "Hopf Algebras," Benjamin, New York, 1969.
\bibitem{TW}
E. Taft, R. Wilson, On antipodes in pointed Hopf algebras, 
{\it J. Algebra} {\bf 29} (1974), 27 -- 32.
\end{thebibliography}
\end{document}